\documentclass{article}

\usepackage{amssymb,theorem}

\title{Sequences and filters of characters characterizing
   subgroups of compact abelian groups}

\date{}
\author{
M. Beiglb\"ock\thanks{mathias.beiglboeck@tuwien.ac.at},
 C. Steineder\thanks{christian.steineder@tuwien.ac.at} \
and R. Winkler\thanks{reinhard.winkler@tuwien.ac.at} \\ \\
{\small Institute of Discrete Mathematics and Geometry,}
\vspace{-1mm}
\\
\vspace{-1mm} {\small TU Vienna}, {\small Wiedner Hauptstra\ss e
8-10,}
\\{\small 1040 Wien, Austria}}

\newtheorem{theorem}{Theorem}[section]

\newtheorem{lemma}[theorem]{Lemma}

\newcommand{\R}{\mathbb{R}}
\newcommand{\T}{\mathbb{T}}

\newcommand{\Z}{\mathbb{Z}}
\newcommand{\N}{\mathbb{N}}
\newcommand{\F}{\mathcal{F}}
\newcommand{\Gd}{\widehat{G}}
\newcommand{\Gdd}{\widehat{G_d}}
\newcommand{\ee}{\varepsilon}
\newcommand{\la}{\alpha_1, \ldots ,\alpha_t}
\newcommand{\Bohr}{\mbox{$B_{(\la,\ee)}$}}

\newcommand{\lga}{\mbox{$\gamma_1, \ldots ,\gamma_d$}}
\newcommand{\lnu}{\mbox{$\nu_1, \ldots ,\nu_l$}}
\newcommand{\sign}{\mbox{sign}}
\newcommand{\g}{\mbox{$\bf{g}$}}
\newcommand{\U}{\mbox{$\mathcal{U}$}}
\newcommand{\lb}{\beta_1, \ldots ,\beta_t}

\newcommand{\seq}[1]{{( {#1}_{n})_{n=1}^\infty}}
\newcommand{\io}[1]{_{#1=1}^\infty}

\newcommand{\limm}[1]{\lim_{{#1}\rightarrow\infty}}

\newcommand{\nhat }[1]{\{1,2,\ldots,#1\}}
\newcommand{\ohat }[1]{\{0,1,\ldots,#1\}}
\renewcommand{\subset}{\subseteq}
\renewcommand{\supset}{\supseteq}
\newcommand{\eps}{\varepsilon}

\begin{document}

\maketitle

\begin{abstract}
Let $H$ be a countable subgroup of the metrizable compact abelian
group $G$ and $f:H\rightarrow \T=\R/\Z$ a (not necessarily continuous) 
character of $H$.   Then there exists a sequence
$\seq \chi$ of (continuous) characters  of $G$ such that 
$\limm n \chi_n(\alpha) = f(\alpha) $ for all $\alpha\in H$ 
and $(\chi_n(\alpha)) \io n$ does not converge whenever
$\alpha\in G\setminus H$. If one drops the countability and 
metrizability requirement one can obtain similar results by 
using filters of characters instead of sequences.
  Furthermore the introduced methods allow to answer questions of
    Dikranjan et al.\\
    \\
    {\footnotesize
    }
\end{abstract}

\noindent {\it Acknowledgements:} The authors are very grateful to
Andr\'as Bir\'o and Vera T. S\'os for substantial
discussions concerning this article. Moreover they thank
Dikran Dikranjan and Gabriel Maresch for valuable remarks
and the Austrian Science Foundation FWF for the
support through project no. S8312.

\section{Introduction}

\subsection{Motivation}

In \cite{biro1} several techniques have been developed to prove
the existence of sequences $(k_n)\io n$ of positive integers
characterizing countable subgroups $H$ of the circle group $\T =
\R/\Z$ in the sense that for  $\alpha \in \T$,
$$\alpha \in H\ \Longleftrightarrow\
   \lim_{n\rightarrow \infty}  k_n \alpha  = 0.$$
These methods were extended in \cite{biro2} to show that if
 $H$ is generated freely by finitely many elements,
a characterization is
possible in an even stronger sense: One can choose a
characterizing sequence such that
$\sum_{n=1}^\infty \| k_n \alpha\| < \infty$
for $\alpha \in H$, while
$\limsup_{n\rightarrow \infty} \|k_n \alpha\| \geq 1/4$ for
$\alpha\in \T \setminus H$. (For $x = r +\Z \in \R/\Z$, $r \in \R$,
 the norm $\| x \|$
denotes the distance between $r$ and the nearest integer.)

In \cite{win} arbitrary subgroups of $\T$ were characterized by filters
on its dual $\Z$. This approach was used in \cite{beig}
to extend the results from \cite{biro2}.

A different approach to
the characterization of finitely generated dense subgroups of
compact abelian groups by sums has recently been introduced in
\cite{biro3} and \cite{biro4}.

Dikranjan et.\ al.\ investigated related questions concerning the
characterization of subgroups of more general topological abelian
groups $G$ (cf. \cite{bar1}, \cite{bar2}, \cite{dik1},
\cite{dik2}). In this article we lift the techniques of
\cite{beig} to this general setting and answer questions stated
in \cite{bar2} and \cite{dik2}.

Our Theorems \ref{char} and \ref{gclosed} have contemporaneously
and independently been proved by Dikranjan and Kunen, cf.\
\cite{dik4}. They also treat questions related to descriptive set
theory. More results in this directions were, for instance, also
obtained by Elia\v{s} in \cite{eli}.

\subsection{Content of the paper}

In Section \ref{Sfilter} we modify the filter method from
\cite{win} for our purposes. Theorem \ref{filt1} essentially states
that arbitrary subgroups of compact abelian groups $G$ can be
characterized by filters on the (discrete)
Pontryagin dual $\Gd$ of $G$. (Such filters are intended to be
the neighborhood filters of 0 w.r.t.\ precompact group topologies
on $\Gd$.)

In Section \ref{Schar} this characterization is used to prove
Theorem \ref{char}:
Among the compact abelian groups $G$ exactly
the metrizable ones have the property that every countable subgroup
$H$ is characterizable in the sense that there is a sequence of
characters $\seq \chi$ in $\Gd$ such that
$$\alpha\in H \ \Longleftrightarrow\ \limm n  \chi_n(\alpha) =0.$$
 This solves Problem 5.3 from \cite{dik2}.
Theorem \ref{gclosed} (which solves Problem 5.1 and Question 5.2
from \cite{dik2}) states that in an arbitrary compact abelian
group every countable subgroup is the intersection of subgroups
characterizable in the above sense.

Section \ref{Sthick} is motivated by Question 5.2 from
\cite{bar2}: Consider the case $G=\T$ and $\Gd = \Z$. It was
established in
 \cite{biro1} that for any countable $H \le \T$ there exists a
sequence of integers $k_1 < k_2< \ldots$ characterizing $H$. Is
it possible to choose the $k_n$ in such a way that the quotients
$\frac{k_{n+1}}{k_n}$ are bounded? We answer this question
affirmatively by proving a stronger assertion, cf.\ Theorem
\ref{thick}. It states that, in a certain sense, characterizing
sequences can be arbitrarily close to having positive density.
This seems to be best possible insofar as (apart from trivial
cases) characterizing sequences always have density 0. Theorem
\ref{thin} is a counterpart to Theorem \ref{thick} and describes
 how sparse characterizing sequences of
 subgroups of $\T$ can be.

In  Section \ref{Sappl}
 we introduce a
refined characterization of subgroups of a compact
metrizable group $G$ by
sequences: For a sequence $\seq \chi$ in $\Gd$ we consider
the set $H$ of
all $\alpha \in H$ for which  $(\chi_n(\alpha)) \io n$
 converges (not necessarily to
$0 \in \T$). $H$ is easily seen to be a subgroup
of $G$ and the
pointwise limit is a (not necessarily continuous) homomorphism
$f: H \to \T$. Theorem \ref{appl} gives a complete
description of the situation: Given any subgroup $H$ of a
metrizable compact abelian group $G$ and any homomorphism $f: H
\to \T$ there is a sequence $\seq \chi$ in $\Gd$ such that
$\chi_n \to f$ pointwise on $H$. If $H$ is countable then one can even
achieve that $H$ is exactly the set of convergence.
If $G$ is a compact (not necessarily metrizable) group and
$H$ is an arbitrary (not necessarily countable) subgroup of $G$,
this result is still valid  when
the convergence of sequences is replaced
by the more general convergence of  filters.
By considering the trivial homomorphism $f \equiv 0$ we see that
Theorem \ref{appl} nicely extends Theorem \ref{char}.
Furthermore this result
allows to construct  counterexamples to Question 5.4 from
\cite{bar2}. (For the complete statement of this question cf.\
Section \ref{Sappl}.)

\subsection{Conventions and notation}

 If not stated
otherwise,
 $G$ is always an infinite locally compact
abelian group.  (For finite $G$ most assertions turn out to be
trivial.)
 Since we are only interested
in abelian groups,
 we use additive
notation. In particular we shall
do so in the group $\T = \R/\Z$ where characters $\chi:
G \to \T$ take their values. Elements of $G$ will be denoted by
$\alpha,\beta,\ldots$. If $H$ is a (not necessarily closed)
subgroup of $G$, we write $H \leq G$. If $A \subset G$ is any
subset, $\langle A \rangle$ denotes the subgroup generated by
$A$. For finite  $A = \{\alpha_1, \ldots, \alpha_n\} \subset G$
 and $M \in
\N$ we put $\langle A \rangle_M:=
  \{ \sum_{i=1}^n k_i \alpha_i : k_i \in \Z,
|k_i| \leq M\}$.

Recall that a group compactification of (any
topological group) G is defined to be a pair $(\iota,C)$ where
$C$ is a compact group and
$\iota: G \mapsto C$ is a continuous homomorphism with
dense image.
 Relative topologies on
$G$ induced by group compactifications are called precompact.
The so called Bohr
compactification $(\iota_{bG}, bG)$ of $G$ is the
compactification of $G$ which is  maximal
in the sense that for
each compactification $(\iota,C)$ of $G$ there is a continuous
homomorphism $\phi: bG \rightarrow C$ with $\phi \circ \iota_{bG} =
\iota_C$.

 $\Gd = \{ \chi : \chi \mbox{ is a
continuous homomorphism from } G \mbox{ to } \T\}$
 denotes the dual group of $G$,
equipped with the compact open topology.
By Pontryagin's Duality
Theorem (cf. for instance \cite{dik3} or \cite{hew}) we know
that for LCA groups $G \cong \widehat{\widehat{G}}$ in the
algebraic as well as in the topological sense via the canonical
mapping $\alpha \mapsto x_{\alpha}$, $x_{\alpha}: \chi \mapsto
\chi(\alpha)$.

We take  $G_d$ to be $G$ endowed with the discrete topology. Duality
theory can be applied to construct the  Bohr
compactification $bG$ of $G$ by setting $bG := \widehat{(\Gd)_d}$
and $\iota_B: \alpha\mapsto x_\alpha$. Accordingly, the
Bohr compactification of $\Gd$ is $\Gdd$. It is natural to call
the precompact topology on $G$ induced by $bG$ the Bohr topology.
On the dual group $\Gd$ the Bohr  topology can be described by the
so called Bohr sets which are defined by
\[
    \Bohr := \left\{ \chi \in \Gd \: : \: \| \chi(\alpha_i) \|
    \leq \ee \mbox{ for } i\in \nhat t \right\},
\] where $\alpha_1,\ldots, \alpha_t \in G$
and $\ee>0$.
These  sets generate the
neighborhood filter  of $0$ in $\Gd$ endowed with the Bohr topology.
Further we put $\Bohr(E) := \Bohr \cap E$ for $E
\subset \Gd$. For $\alpha \in G$ and $B \subset \Gd$ we write
$\|\alpha B \| := \sup\{ \| \chi(\alpha) \| \: : \: \chi \in
B\}$.

\section{Characterizing filters}\label{Sfilter}

We will make use of so called filter limits: Let $S$
be any set, let $\F$ be a filter on $S$, let $y$ be a point in a
topological space $X$ and let $f:S\rightarrow X$ be function.
Then $$\F-\lim_s f(s)=y $$ iff for every neighboorhood $U$
of $y$, $\{s\in S :f(s)\in U\}\in \F$.

We remark that  filter limits are more general then
limits along sequences: For a sequence
$\seq x$ in $X$ put
$$\F_{\seq x} =\{A\subset X:\exists\ m\in \N \ \mbox{such that}\
    \{x_n:n\geq m\} \subset A\} $$
Then  $\F_{\seq x}-\lim_s f(s)$ exists iff  
$\limm nf(x_n)$ exists and in this case they coincide. \\

Let $H \leq G$ be a  subgroup of the compact abelian group $G$.
Our task is to show that $H$ can be characterized by a filter
$\F_H$  on $\Gd$ in the sense that we have
 $\F_H-\lim_\chi \chi(\beta)=0$ iff $\beta\in H$.
It is clear that for all $\alpha\in G$ and all $\ee>0$
the set $B_{(\alpha,\ee)}$ has to be an element of   $\F_H$
to assure convergence for elements of $H$. 
By the filter properties
of $\F_H$ the intersection of finitely many
such sets will again be an element
of $\F_H$. 
Thus it would be natural to define $\F_H$ to be the filter generated
by the sets $\Bohr$ where $\la \in H, \eps>0$. 
This definition yields the minimal filter with the required property
and corresponds to the precompact group topology on $\Gd$ induced 
by $H$. 
Later it will be important to us that we may neglect finite sets
of characters. Therefore we will also take all cofinite sets to
be elements of $\F_H$. This leads to the following definition:
 \[
    \F_H := \left\{
        F \subseteq \Gd  :
        \begin{array}{l}
         \exists \: \la \in H, \ee > 0, \Gamma\subset \Gd,
               |\Gamma|<\infty  \\
         \mbox{such that } \Bohr(\Gd \setminus \Gamma) \subseteq F
         \end{array}
    \right\}.
\]

\begin{theorem} \label{filt1}
    Let $G$ be an infinite  compact abelian group,
let $H$ be a subgroup of $G$ and let the filter $\F_H$ be defined as above.
Then for all $\beta\in G$
    \[
        \F_H - \lim_{\chi}  \chi (\beta)  = 0 \:
    \Longleftrightarrow  \:  \beta \in H.
    \]
\end{theorem}
In the course of the proof we will employ the following lemma which
will also be useful later on:

\begin{lemma}\label{lem0}
Let $G$ be a compact abelian group. Then $\Gd$ is dense in $\Gdd$
w.r.t.\ pointwise convergence. Thus, for any countable subset $H$ of
$G$ and any $\chi \in \Gdd$ there exists a sequence $\seq \chi$
 in $\Gd$
such that $\chi_n(\alpha) \to \chi(\alpha)$ $(n \to \infty)$
for all $\alpha \in H$.
\end{lemma}
{\it Proof:} As explained in the introduction, the compact group
$\Gdd$ is, with the set theoretic inclusion as dense embedding,
the Bohr compactification of the discrete group $\Gd$. This
proves the first part. Thus for $H = \{\alpha_1,\alpha_2,\ldots\}
\subset G$ and every $n \in \N$ there is a $\chi_n \in \Gd$ with
$\| \chi_n(\alpha_i) - \chi(\alpha_i)\| < \frac 1n$ for all
$i\in \nhat n$. It follows that $\chi_n \to \chi$ pointwise on
$H$. $\Box$
\\

{\it Proof of Theorem \ref{filt1}:}
    The definition of $\F_H$
     guarantees that
    $\F - \lim_{\chi} \| \chi (\beta) \| = 0$ for all
    $\beta \in H$. For the converse we prove that, given $\beta \not\in
    H$, for all $\la \in H$ and every $\ee > 0$ there exist
    infintely many characters
    $\chi \in \Bohr$ with $\| \chi(\beta) \| \ge 1/4$ which implies
    that $\{\chi \in
    \Gd : \| \chi (\beta) \| < 1/4 \} \not\in \F$.
    First we see that there exists at least one such character:
    Consider the Bohr compactification $\Gdd$ of $\Gd$. $\Gdd$ separates
    subgroups and points of $G$. Hence there exists some
    $\phi\in \Gdd$ such that
    \[
        \phi(\alpha) = 0 \mbox{ for all } \alpha \in \langle \la
        \rangle \quad \mbox{and} \quad c= \phi(\beta) \not= 0,
    \]
    w.l.o.g. $\| \phi(\beta) \| \ge 1/3$ (otherwise take
    an appropriate multiple $2\phi,3 \phi,\ldots$).
    By Lemma \ref{lem0} $\phi$ can be approximated arbitrarily
    well on finitely many points by a character.
    Thus we find some   $\chi \in \Gd$
    such that $\|\chi(\alpha_i) \|\leq \ee, 1\leq i \leq t,
    \| \chi(\beta)\| >1/4$.

    Next we prove that for $\eps>0$ each $\Bohr$ contains infinitely
    many  $\chi$ with $\| \chi(\beta)\| \ge 1/4$.  Let
    $U := \{ (\chi(\alpha_1), \ldots,
        \chi(\alpha_t), \chi(\beta)) :
          \chi \in \Gd \} \leq \T^{t+1}$.
    We distinguish two cases: \\
    1. $U$ is finite, say $U = \{ u_1, \ldots, u_k \}$.
    There is some $i$, say $i=1$, with $u_1 := (0,\ldots,0,c)$.
    Then the sets
    \[
        \Upsilon_i := \left\{ \chi \in \Gd :
        (\chi(\alpha_1), \ldots, \chi(\alpha_t),\chi(\beta)) = u_i
        \right\}, \quad i= 1,\ldots,k,
    \]
    and particularly $\Upsilon_1$ are infinite, or\\
    2. $U$ is an infinite subgroup of $\T^{t+1}$. But then each point of
    $U$ is an accumulation point.\\
    In both cases we find infinitely many $\chi$ with the required property.
    $\Box$

\section{Characterizing  countable subgroups} \label{Schar}

We solve Problem 5.3 from \cite{dik2}: For which
compact abelian $G$ can every countable subgroup $H$ be
characterized by a sequence of characters?

For $A \subset \Gd$ we write $\lim_{\chi \in A} \chi(\beta)
=0$  iff $\{ \chi \in A :
\chi(\beta) \ge \ee\}$ is finite for all $\ee >0$. (I.e. instead of
the characterizing sequence $\seq \chi$ we consider the characterizing
 set
$A=\{\chi_n:n\in \N\}$.)

\begin{theorem} \label{char}
    Let $G$ be an infinite compact abelian  group and let  $H \leq G$ be a
    countable subgroup. Then the following statements are
    equivalent:
\begin{itemize}
    \item[(i)] $G$ is metrizable.
    \item[(ii)]  There exists a countable set $A \subset \Gd$, such that
    \[
        \beta \in H \: \Longleftrightarrow \: \lim_{\chi \in A}
        \chi(\beta) =0 .
    \]
\end{itemize}
\end{theorem}

{\bf Remark.} The proof of Theorem \ref{char} actually
 shows that (if $G$ is metrizable) for every $\sigma < 1/3$
the characterizing set $A$ can be chosen in such a way that $\beta
\not \in H$ implies $\limsup_{\chi \in A} \|\chi(\beta)\| \ge \sigma$.
Using a diagonalization argument it is not difficult
to achieve  $\limsup_{\chi \in A} \|\chi(\beta)\| \ge 1/3$
 and it is easy to see that this is best possible.
\\

The proof of $(i) \Longrightarrow (ii)$ employs several lemmas
which we formulate now and verify at the end of this section.
According to our assumptions, in these lemmas $G$ is an infinite compact
abelian metrizable group.

\begin{lemma} \label{lem1}
    Let $\tau \in \T$ and $n \in \N$.
 Assume that $\| i \tau \| \leq \sigma < 1/3$
     for all $i\in \nhat n$. Then $\| \tau \| \leq \sigma/n$.
\end{lemma}

\begin{lemma} \label{lem11}
    Assume that $\gamma_1, \ldots ,\gamma_d \in G$ freely generate
    $H \leq G$. For arbitrary nonempty open sets $I_1, \ldots, I_d$ in
    $\T$ there exists $\chi \in \Gd$ such that $ \chi(\gamma_i )  \in
    I_i$ for all $i\in \nhat d$.
\end{lemma}

\begin{lemma}       \label{lem4}
    Let  $\la \in G$, $\ee >0$ and $\sigma < 1/3$.
    \begin{enumerate}
    \item
    For all finite $\Gamma \subset \Gd$ and all $\beta \in G$
        \[
            \| \beta  \Bohr(\Gd \setminus \Gamma) \| \leq \sigma
            \ \Longrightarrow\ \beta \in
            \langle\la\rangle.
        \]

    \item Moreover there exists $M\in \N$ such that
    for all finite $\Gamma \subset \Gd$ and all $\beta \in G$
        \[
            \| \beta  \Bohr(\Gd \setminus \Gamma) \| \leq \sigma
            \ \Longrightarrow\ \beta \in
            \langle\la\rangle_M.
        \]

    \item If  $V\supset \langle \la\rangle_M$ is an open subset of
          $G$ then for all
    finite $\Gamma \subset \Gd$
    there exists a finite set $E\subset \Gd\setminus \Gamma$ such
    that for  $\beta \in G$
        \[
            \| \beta  \Bohr(E) \| \leq \sigma
            \ \Longrightarrow\ \beta \in
            V.
        \]
    \end{enumerate}
\end{lemma}

\begin{lemma}   \label{lem3}
    Let $R_1 \subseteq R_2 \subseteq \ldots$ be finite subsets
    of $G$. There exists a sequence of open sets $V_n \subset G, n\in\N$
    such that
    \begin{enumerate}
        \item
            $V_n \supseteq  R_n$.
        \item
            $ \lim\inf_{n\rightarrow\infty} V_n
           = \bigcup_{m = 1}^\infty \bigcap_{n = m}^\infty
            V_n = \bigcup_{n = 1}^\infty R_n$.
    \end{enumerate}
\end{lemma}

\noindent {\it Proof of Theorem \ref{char}:}

$(i) \Longrightarrow (ii)$:
We will first construct the  set $A \subset \Gd$ and then
 prove  that
$\beta \in H$ iff $\lim_{\chi \in A} \chi(\beta) =0$.

 Let $H =: \{ \alpha_t : t \in \N \}$
and pick $\ee=\sigma\in (0,1/3)$.  Using \ Lemma \ref{lem4},2
we can choose a sequence $(M_t)_{t=1}^\infty$ such that for every
finite $\Gamma\subset \Gd$ and all $\beta \in G$
        \[
            \| \beta  \Bohr(\Gd \setminus \Gamma) \| \leq \ee
            \Longrightarrow \beta \in
            \langle\la\rangle_{M_t}.
        \]
Next put for $t\in \N, R_t := \langle \la \rangle_{M_t}$ and
define $V_t \supseteq R_t$ according to Lemma \ref{lem3}
such that $\liminf_{t\rightarrow \infty} V_t = H$.

Using Lemma
\ref{lem4},3 we choose a  finite set $E_1 \subset \Gd$
 such that $ \| \beta
B_{(\alpha_1,\ee)}(E_1) \| \leq \ee$ implies $\beta \in V_1$.
By employing Lemma \ref{lem4},3 again, we find
$E_2 \subset \Gd\setminus E_1$ such that
$\|\beta
B_{(\alpha_1,\alpha_2,\ee)}(E_2) \| \leq \ee$ implies $\beta \in V_2$.
Continuing in this fashion we arrive at a sequence $(E_t)_{t=1}^\infty $
of disjoint subsets of $\Gd$ such that for each $t\in \N$
$$\| \beta B_{(\alpha_1,\ldots,\alpha_t,\ee)}(E_t) \| \leq \ee
\ \Longrightarrow \ \beta \in V_t.$$
Finally we put
$
    A:= \bigcup_{t = 1}^\infty \Bohr(E_t).
$

Assume that $\beta \in H$. To prove $\lim_{\chi \in A} \|\chi(\beta)\|
=0$ note that, for arbitrary $n \in \N$, there exists $T=T(n) \in
\N$  such that $i \beta \in \{ \alpha_t : t \leq T\}$ for
all $i\in \nhat n$. Thus whenever $\chi \in \Bohr(E_t)$
for some $t\geq T$
we have $\| i \chi(\beta) \| \leq \ee$ for all $1 \leq i \leq n$.
By Lemma \ref{lem1} this yields $\| \chi(\beta) \| \leq \ee/n$.
Since $n$ was arbitrary we get $\lim_{\chi \in A}
\|\chi(\beta)\| =0$.

Conversely assume that $\limsup_{\chi \in A} \|
\chi (\beta) \| < \ee$ for some $\beta \in G$.
Then  for all but finitely
many $t\in \N$ we have $\|\beta \Bohr (E_t)\| \leq \eps $.
Thus there exists $t_0 \in \N$ such that
$\beta \in V_t$ for all $t \ge t_0$ which yields $\beta \in H$ by the
choice of the sequence $(V_t)\io t$.\\

$(ii) \Longrightarrow (i)$: Let $H \leq G$ be an arbitrary
countable subgroup characterized by the countable set $A 
 \subset \Gd$.  Define $\Lambda := \langle A
\rangle$ and
\[
    \Lambda^0  := \{ g \in G : \chi(g) = 0
    \mbox{ for all } \chi \in \Lambda \},
\]
the annihilator of $\Lambda$. Clearly $\Lambda^0 \leq H$, thus $|
\Lambda^0 | \leq \aleph_0$. Since
$\widehat{\Lambda} \cong G/ \Lambda^0$
we have $ w( G/\Lambda^0) = w(\widehat{\Lambda}) = |\Lambda| = \aleph_0$,
where $w$ denotes the topological weight, i.e. the least cardinal number
of an open basis (cf.\ \cite{hew}, 24.10 and 14.14). Hence $G/\Lambda^0$ and
$\Lambda^0$ have at most countable weight and therefore also $G$
(\cite{hew}, 5.38), implying that $G$ is metrizable. $\Box$\\

Let $G$ be a compact abelian group.
In \cite{dik2} subgroups  characterized by a sequence
$\seq \chi$ in $\Gd$ are denoted by
\[
    s_{\seq \chi} (G) := \left\{ \alpha \in G : \limm n \chi_n (\alpha)=
    0 \right\}.
\]
Furthermore such subgroups are called basic
$\g$-closed subgroups.
 According to Theorem \ref{char}
every countable subgroup of $G$ is basic $\g$-closed iff $G$ is
metrizable.

A group $H \leq G$ is called $\g$-closed if
it is representable as the intersection of basic
$\g$-closed subgroups.
The next theorem deals with $\g$-closed subgroups and solves
Problem 5.1 from \cite{dik2}.

\begin{theorem} \label{gclosed}
    Every countable subgroup $H$ of a compact abelian group $G$
    is $\g$-closed.
\end{theorem}
{\it Proof:} For arbitrary $\beta \in G \setminus H$ there is a
$\chi \in \Gdd$ with $\chi(\alpha)=0$ for all $\alpha \in H$ and
$\| \chi(\beta)\| \ge \frac 13$. Thus Lemma \ref{lem0}
immediately yields a sequence of  $(\chi_n^{\beta}) \io n$ in
$\Gd$ characterizing a subgroup
\[
    H_{\beta}:= s_{({\chi_n^{\beta}})\io n} (G) =
 \left\{ \alpha \in G :  \limm n  \chi_n^{\beta}
(\alpha)  =0 \right\} \leq G
\]
with $H \le H_{\beta}$ and $\beta \not \in H_{\beta}$.
Thus $H = \bigcap_{\beta \in G \setminus H} H_{\beta}$. $\Box$

\paragraph{Proofs of Lemmas \ref{lem1} to  \ref{lem3} :\\\\}

\noindent We assume the group $G$ to be compact abelian and
metrizable. Lemma \ref{lem1} is elementary, so we  skip the proof.\\
\\
{\it Proof of Lemma \ref{lem11}:} Assume that
\[
    A := \Gd \langle\gamma_1, \ldots,
    \gamma_d\rangle = \{ \chi(\alpha) \: : \: \chi \in \Gd, \alpha
    \in \langle\gamma_1, \ldots, \gamma_d\rangle \}
\]
is not dense in $\T^d$, i.e.\ $H := \overline{A} \lneq \T^d$.
There is a nontrivial character of $\T^d$ vanishing on $H$, i.e.\
a nonzero vector $h=(h_1,\ldots,h_d) \in \Z^d$
such that $\sum_{i=1}^d h_i x_i =0$ holds for all
$x=(x_1, \ldots , x_d) \in H$. Fix an arbitrary $\chi \in \Gd$
and put $x_i = \chi(\gamma_i)$. Then
\[
    0 = \sum_{i=1}^d h_i \chi( \gamma_i) =  \chi\left(\sum_{i=1}^d h_i
    \gamma_i\right).
\]
Since this holds for all $\chi \in \Gd$ we have
$\sum_{i=1}^d h_i \gamma_i = 0$, contradicting the independence
of the free generators $\gamma_i$, $1 \leq i \leq d$.
$\Box$\\
\\

{\it Proof of Lemma \ref{lem4}:}
Let  $B_0:= \Bohr(\Gd \setminus \Gamma)$.

1. Let $\F=\F_{\langle \la \rangle}$ be the filter of Theorem
\ref{filt1} characterizing $\langle \la \rangle$ and let $\delta
>0$ be arbitrary.  Under the assumption
$\| \beta  B_0\|\leq \sigma < 1/3$ we
have to show that
\[
    F_{\delta} := \left\{ \chi \in \Gd \setminus \Gamma
     : \| \chi (\beta) \| \leq \delta
    \right\} \in \F.
\]
Choose $m \in \N$ such that $\delta \ge \sigma /m$ and
let $B_1:= B_{(\la,\ee/m)}(\Gd
\setminus \Gamma)$. By definition of $\F$ we have
$B_1\in\F$. For all $\chi \in B_1,$ $i\in \nhat m$, we have
$i\chi \in B_0$. Thus  $\|i \chi
(\beta)\| \leq \sigma$ for all $i \in \nhat m$ and Lemma \ref{lem1}
yields $\| \chi (\beta) \| \leq \sigma/m < \delta$. Thus
$B_1 \subset F_{\delta}$ and hence $F_{\delta} \in \F$.
\\

2. Assume that $H:= \langle \la \rangle$ is infinite (otherwise
the assertion follows immediately). Since $H$ is a finitely
generated abelian group there exists a decomposition $H = T
\oplus F$ where $F$ is freely generated by $\lga$  and $T=\langle
\lnu \rangle = \oplus_{i=1}^l \langle \nu_i \rangle$ is the
torsion subgroup of $H$. Hence $\langle \nu_i \rangle \cong
\Z_{e_i}$ for some $e_i \in \N$ and
\[
    \langle \la \rangle = \langle \lga \rangle \oplus \langle \lnu \rangle
    \cong \Z^d \oplus \textstyle{\bigoplus_{i=1}^l} \Z_{e_i}.
\]
Let $\delta >0$ be such that $\| \chi(\gamma_i) \| \leq \delta$
for $i \in \nhat d$ and $\chi  (\nu_j) =0 $ for  $j
\in \nhat l$ implies $\| \chi (\alpha_k)\| \leq \ee$ for
$k \in \nhat t$. Pick now any $\beta \in G$ with $\|
\beta B_0 \| \leq \sigma < 1/3$. By 1. above
 we have $\beta \in H$, thus $\beta =
\sum_{i=1}^d r_i \gamma_i + \sum_{j=1}^l s_j \nu_j$ for some
$r_i
\in \Z$, $s_j\in \ohat {e_j-1}$, $i\in \nhat d,j\in \nhat l$.
 Let $e := \prod_{j=1}^l e_j$. By Lemma \ref{lem11} there
exist infinitely many $\chi \in e \Gd := \{ e \chi' : \chi'
\in \Gd \}$ such that
\[
    \sign(r_i) \chi(\gamma_i) \in \left[ \frac{1}{3 \sum_{j=1}^d |r_j| },
     \frac{2}{3 \sum_{j=1}^d |r_j| } \right] + \Z
\]
holds for $i\in \nhat d$. Therefore we have
\[
   \ r_i \chi(\gamma_i)  \in \left[\frac{|r_i|}{3 \sum_{j=1}^d |r_j|} ,
    \frac{2 |r_i|}{3 \sum_{j=1}^d |r_j| } \right] + \Z
\]
 for all $i\in \nhat d$. Summing up and using that
$\chi(\nu_j) =0$ for $j\in \nhat l$ this leads to
\[
     \chi(\beta) = \sum_ {i=1}^d \ r_i \chi(\gamma_i)
     \in \left[ \frac{1}{3},  \frac{2}{3} \right] + \Z.
\]
Thus $\chi \not \in \Bohr$ and hence there is $j \in
\{1,\ldots,t \}$ with $\| \chi (\alpha_j) \| > \ee$ and therefore
$\delta < \frac{2}{3 \sum_{j=1}^d |r_j|}$. Equivalently
$\sum_{i=1}^d |r_i| < \frac{2}{3 \delta}$.
 So there are only finitely many choices for $\beta$ and we may
put an
universal bound $M$ on the coefficients in the linear combination
  $\beta =\sum_{i=1}^r k_i \alpha_i$.
\\

3. Clearly, the set
\[
    I:=\{\beta \in G:\| \beta B_0\|\leq \sigma\}= 
    \bigcap_{\chi \in B_0}\{ \gamma \in G  :
    \|\chi (\gamma)\| \leq \sigma\}
\]
is closed and by 2. we have $I\subset \langle \la\rangle_M\subset V$.
Thus $I\cap V^c=\emptyset.$
 By compactness of $G$ there exists a
finite set $E \subset B_0$ such that
\[
    \bigcap_{\chi \in E}\{ \gamma \in G :
    \|\chi (\gamma)\| \leq \sigma\} \cap V^c = \emptyset.
\]
This $E$ is as required. $\Box$
\\

{\it Proof of Lemma \ref{lem3}:} Let $\rho$ be a metric on $G$
compatible with its topology. Since the sets $R_1 \subseteq R_2
\subseteq \ldots \subseteq G$ are finite there is a sequence
$\seq d$ of positive reals decreasing to $0$ such that
\begin{eqnarray*}
    2 d_n &<& \min\{ \rho(\alpha,\alpha') \: : \: \alpha,\alpha' \in
    R_n, \alpha \not= \alpha'\}\\
    d_n + d_{n+1} &<& \min \{ \rho(\alpha, \alpha')\: : \: \alpha
    \in R_n, \: \alpha' \in R_{n+1} \setminus R_n \}.
\end{eqnarray*}
Define
\[
    V_n := \{ \beta \in G \: : \: \exists \: \alpha \in R_n \mbox{ with }
    \rho(\beta, \alpha) < d_n \}.
\]
By monotonicity of the sets $R_n$, $\beta \in \bigcup_{n= 1}^\infty
R_n$ implies $\beta \in \bigcup_{m = 1}^\infty \bigcap_{n = m}^\infty V_n$.
Conversely, assume $\beta \in \bigcup_{m = 1}^\infty \bigcap_{n= m}^\infty
V_n$ or, equivalently, that there exists an $m$ with $\beta \in
V_n$ for all $n \ge m$. According to the definition of the sets
$V_n$ there exists a unique $\alpha_n \in R_n$ such that
$\rho(\alpha_n,\beta) < d_n$ for $n \ge m$. Moreover the choice
of the $d_n$ guarantees that $\alpha_m = \alpha_{m+1}  = \ldots$
and so $\rho(\beta,\alpha_m) = \rho(\beta,\alpha_n) \leq d_n \to
0$. Hence $\beta = \alpha_m \in R_m \subseteq \bigcup\io n R_n$.
$\Box$\\

\section{Thick and thin characterizing sequences}\label{Sthick}

Question  5.2 from \cite{bar2} asks: Does every
countable subgroup $H$ of $\T$ admit a characterizing sequence
$\seq k$ with bounded quotients, i.e.\ $q_n =
\frac{k_{n+1}}{k_n} \le C$ for all $n \in \N$ and some $C \in
\R$?

 We  answer this question affirmatively by
proving a stronger result. Which type
of statement can be expected?
Assume that $\alpha \in H$ is irrational. Then, by
uniform distribution of the sequence $({n\alpha})\io n$, the set of
all $k \in \N$ with $\| k\alpha \| < \ee$ has density $2 \ee$. Thus
(with the exception of trivial cases) characterizing sequences
have zero density. Furthermore the length of their gaps is tends to 
infinity. In particular the thickest characterizing sequences we can expect
might have a density which converges to zero very slowly in some
sense. This is the content of the following result.

\begin{theorem} \label{thick}
    Let $H \leq \T$ be a countable subgroup and let $(\ee_j)_{j=1}^\infty$
 be a
    sequence with $0 \leq \ee_j \leq 1$ that converges to $0$
    (arbitrarily slowly). Let $\N$ be partitioned into nonempty intervals
    $I_j = \{i_j, i_j+1, \ldots ,i_{j+1}-1\}$ with $i_0 =0$ and
    $\lim_{j\rightarrow \infty} ( i_{j+1}- i_j) = \infty$.
Then there exists a
    sequence $\seq k$ of nonnegative integers characterizing $H$
    such that
    \[
        \frac{ \left| \{n: k_n \in I_j\} \right|}{| I_j|} \ge \ee_j
        \quad \mbox{ for all } j.
    \]
\end{theorem}
{\it Proof:} Let, according to Theorem \ref{char} (or to
\cite{biro1}), $c_1 < c_2 < \ldots \in \N$ be any sequence
characterizing $H$. We are going to construct a sequence $d_1 <
d_2 < \ldots \in \N$ containing at least $\ee_j |I_j|$ elements
in each $I_j$ such that $\| k_n \alpha \| \to 0$ for all $\alpha
\in H$. Then $A = \{k_1 < k_2 < \ldots \} = \{d_n : n \in \N \}
\cup \{ c_n : n \in \N \}$ clearly has the desired properties.

Let $H =: \{ \alpha_t  :  t \in \N \}$, put $I_j^+ (0) := I_j$ and
$$
    I_j^+(t) := \left\{ k \in I_j : \| k \alpha_i \| <
    1/t \mbox{ for all } i \in \nhat t \right\}
$$
for  $t \geq 1$. For each $j \in \N$
let $t_j$ be the maximal $t \in \ohat j$ such that $|I_j^+(t)| \ge
\ee_j | I_j |$ and put
\[
    \{d_1 < d_2 < \ldots \} =\textstyle{ \bigcup_{j = 0}^\infty}\
    I_j^+ (t_j).
\]
It suffices to show that $t_j \to \infty$ for $j \to \infty$ or,
equivalently, that for each $t_0 \in \N$ there exists  $j_0$
such that for all $j \ge j_0$
\[
    \Big| \left\{ d \in I_j : \| d \alpha_i \| < 1/t_0
    \mbox{ for all } i = 1, \ldots, t_0 \right\} \Big| \ge \ee_j
    I_j.
\]

Since $\ee_j \to 0$ for $j \to \infty$ this is an immediate
consequence of the well distribution (cf. \cite{kuip}) of the
sequence $(ng)\io n$ in the closed subgroup $G \leq
\T^{t_0}$ generated by $g=(\alpha_1, \ldots, \alpha_{t_0}) \in
\T^{t_0}$: The open subset $O \subset G$ of all $(\beta_1,\ldots,
\beta_{t_0})$ with $\| \beta_i \| < 1/t_0$ has positive Haar
measure $\mu (O)$ and the set of all $k \in \Z$ with $k g \in O$
has uniform density $\mu(O) > 0$. $\Box$\\

Theorem \ref{thick} indeed answers the question about quotients:
Take, for instance, $i_j = j^2$ and choose the $\ee_j$
in such a way that 
at least one $k_n$ lies in each $I_j$. Then the
quotients $q_n = \frac{k_{n+1}}{k_n}$ tend to 1. 
This example can be modified in many ways.

It has been proved in \cite{bar2} that $q_n \to \infty$ implies
that the corresponding characterized group $H$ is uncountable.
Thus, for a given countable $H$, characterizing sequences cannot
be arbitrarily sparse in this sense. Nevertheless we have:

\begin{theorem} \label{thin}
Let $H$ be a countable subgroup of $\T$ and let  $m_1 < m_2 < \ldots$
be an (arbitrarily fast) increasing sequence of positive integers.
Then there is a characterizing sequence $k_1 < k_2 < \ldots$ for
$H$ with $m_n < k_n$ for all $n \in \N$.
\end{theorem}
{\it Proof:} Let $\seq c$ be any  characterizing sequence of $H$.
Put $k_{2 n} := c_{j_n}$ and $k_{2n+1} := c_{j_n}+c_n$ where
$j_n$ is large enough
in the sense  that $k_{2n} > m_{2n}$ and  $k_{2n+1} >
m_{2n+1}$. Clearly $\alpha \in H$ implies $ k_n \alpha  \to
0$. On the other hand, if $\beta \in \T$ and $k_n \beta \to
0$ then also $ ( k_{2n+1} - k_{2n}) \beta  =  c_n \beta 
\to 0$. $\seq c$ characterizes $H$, therefore $\beta \in H$.
$\Box$\\

Theorem \ref{thin} implies that for any countable $H \leq \T$
there are sequences $\seq k$ characterizing $H$ with
$\limsup_{n\rightarrow \infty}
\frac{k_{n+1}}{k_n}= \infty$:
In Theorem \ref{thin} put $m_n= n^n$  and let $k_1<k_2<\ldots$ be
a characterizing sequence of $H$ such that $m_n \leq k_n$ for 
all $n\in \N$. Then
$$\sup_{n\in \N} \frac{k_{n+1}}{k_n} 
\geq \sup_{n\in \N}\sqrt[n]{\prod_{i=1}^n \frac{k_{i+1}}{k_i}}
\geq \sup_{n\in\N}\sqrt[n]{\frac {k_n} {k_1}}
 \geq \sup_{n\in \N} \sqrt[n]{\frac{n^n}{ k_1}} = \infty.$$ 

For more sophisticated methods to generate sparse characterizing
sequences we refer to \cite{beig} and
\cite{biro2}: E.g.   
for a countable subgroup $H\leq \T$ one can 
construct a characterizing sequence
$\seq k$ such that for all $r>0$ and $\alpha\in H$,
$\sum_{n = 1}^\infty \| k_n \alpha \|^r <
\infty$.\\

The idea of the proof of Theorem \ref{thin} has further remarkable
extensions. We will analyze them more detailed in the next
section.

\section{Groups as sets of convergence}\label{Sappl}

The following theorem presents a generalized approach to the
characterization of subgroups.

\begin{theorem} \label{appl}
    Let $G$ be a compact abelian group.
    \begin{enumerate}
        \item
        Let $\F$ be a filter on $\Gd$. Then  the set $H$ of all
        $\alpha \in G$ for which $\F-\lim_\chi \chi(\alpha)$ exists is
        a subgroup of $G$. The mapping $f: H \mapsto \T$, $\alpha \mapsto
        \F-\lim_\chi \chi(\alpha)$ is a group homomorphism.

        In particular if $\seq \chi$ is a sequence in $\Gd$, the   set $H$ of all
        $\alpha \in G$ for which
        $\lim_{n\rightarrow \infty}\chi_n (\alpha)$ exists is
        a subgroup and the mapping $f: H \mapsto \T$, $\alpha \mapsto
        \lim_{n\rightarrow \infty} \chi_n (\alpha)$
        is a group homomorphism.

        \item
        Let conversely  $H$ be a subgroup of  $G$ and let
        $f: H \mapsto \T$ be a homomorphism. Then
        there exists a filter $\F$ on $\Gd$ such that $\F-\lim_{\chi}
        \chi(\alpha) = f(\alpha)$ for all $\alpha \in H$ and
        $\F-\lim_{\chi} \chi(\beta)$ does not exist whenever
        $\beta \not\in H$.

        \item
        If furthermore $H \leq G$ is countable then there exists a sequence
         $\seq \chi$ in $\Gd$ such that
        \[
            \chi_n (\alpha) \to f(\alpha) \quad \mbox{for all }
            \alpha \in H.
        \]
        \item
        If $G$ is metrizable and $H$ is countable then there exists a
        sequence $(\chi_n') \io n$ in $\Gd$ such that
        \[
            \chi'_n (\alpha) \to f(\alpha) \quad \mbox{for all }
            \alpha \in H
        \]
        and $(\chi'_n(\beta))\io n$ does not converge if $\beta \not \in
        H$.
        \end{enumerate}
\end{theorem}
{\it Proof:} 1: Clear.\\

2. For $\alpha_1,\ldots,\alpha_t \in H$ and  $\ee >0$
put
\[
    F(\la,\ee):=
    \left\{ \chi \in \Gd : \| \chi(\alpha_i) - f(\alpha_i)
    \| \leq \ee \mbox{ for } i = \nhat t \right\}
\]
and
\[
    \F= \F(H,f) := \left\{ F \subset \Gd \: : \:
    \begin{array}{c}
        \exists \:  \la \in H, \: \exists\:  \ee > 0 \\
        \mbox{such that } F(\la,\ee) \subseteq F
    \end{array}
    \right\}.
\]
We have to show that
\begin{itemize}
\item[(a)]   $\F$ is a filter.
\item[(b)]   For all $\alpha \in H$ : $\F-\lim_{\chi}
        \chi(\alpha) = f(\alpha)$.
\item[(c)] For all $\beta \not\in H$ : $\F-\lim_{\chi}
        \chi(\beta)$ does not exist.
\end{itemize}

ad (a): Since the set $F(\la,\lb,\min(\ee_1,\ee_2)) \in \F$ is
contained in $F(\la,\ee_1) \cap F(\lb, \ee_2)$ it suffices to
show that each $F(\la,\ee)$ is not empty.

There exists an extension of $f : H \mapsto \T$ to $\chi: G
\mapsto \T$ such that $\chi \in \Gdd$. By Lemma \ref{lem0} there
is a $\chi' \in \Gd$ such that $\| \chi'(\alpha_i) -
\chi(\alpha_i)\| \leq \ee$ for $i = 1, \ldots,t$. Hence $\chi'
\in F(\la,\ee) \not= \emptyset$.

ad (b): Let $\alpha \in H$ and $U \in \U(f(\alpha))$. There
exists an $\ee > 0$ such that $\| \xi - f(\alpha) \| < \ee$
implies $\xi \in U$. $\chi \in F(\alpha, \ee) \in \F$ implies $\|
\chi (\alpha) - f(\alpha) \| < \ee$ proving $\F-\lim_{\chi}
\chi(\alpha) = f(\alpha)$.

ad (c): Let $\beta \not\in H$ and $F \in \F$ be arbitrary. We will
show that there exist $\chi_1, \chi_2 \in F$ such that $\|
\chi_1(\beta)- \chi_2(\beta)\| \ge 1/4$. $F \in \F$ implies that
there exist $\la \in H$ and $\ee >0$ such that $F(\la,\ee)
\subseteq F$. Note that there is a $\chi' \in \Gdd$ with
$\chi'(h)=0$ for all $h \in H$ and $\chi'(\beta) \ge 1/3$. By
Lemma \ref{lem0} there exists a $\chi \in \Gd$ such that $\|
\chi(\alpha_i)\| < \ee/2$ for $i = 1, \ldots,t$ and $\chi(\beta)
> 1/4$.  Pick $\chi_1 \in F(\la, \ee/2) \subseteq F$
arbitrary and let $\chi_2 = \chi+\chi_1$. Then $\chi_2$ is also in
$F$ and $\| \chi_2(\beta)- \chi_1(\beta) \| = \| \chi(\beta)\| >
1/4$.\\

3. Let $H=\{\alpha_t,t\in \N \}$. The proof of $2.$ shows that
for each $n \in \N$ there is a $\chi_n \in \Gd$ such that $\|
\chi_n(\alpha_i)-f(\alpha_i) \| < 1/n$ for $i= 1, \ldots,n$. The
sequence $\seq \chi$ has the desired properties. \\

4. If $G$ is metrizable and $H$ is countable we know by Theorem
\ref{char} that there exists a sequence
$(\widetilde{\chi_n})_{n=1}^\infty$
in $\Gd$ such that
\[
    \| \widetilde{\chi_n} (\alpha) \| \to 0 \quad \mbox{iff} \quad
    \alpha \in H.
\]
Let furthermore $\seq \chi$ be as in $3.$ and define $\chi'_{2n}
:= \chi_n$ and $\chi'_{2n+1} := \chi_n + \widetilde{\chi_n}$. Then
$\chi'_n (\alpha) \to f(\alpha)$ for all $\alpha \in H$.

Conversely, for $\beta \notin H$ the sequence $(\chi'_n(\beta))\io n$ cannot
converge: If  $\chi'_n (\beta) \to c$ for some $c \in \T$, then
$\widetilde{\chi_n} (\beta) = \chi'_{2n+1}(\beta)-\chi'_{2n}(\beta) \to 0$.
Hence $\beta \in H$, contradiction. $\Box$\\

We want to apply Theorem \ref{appl} to Question 5.4 in \cite{bar2}
which, in our notation, reads as follows.
Let $\seq c$ be a sequence in $\Z$.
Are the subsequent conditions (i) and (ii) equivalent?
\begin{itemize}
\item [(i)] There exists  a precompact abelian group $G \supset\Z$
such that $c_n \to h$ in $G$ and $\langle h \rangle \cap \Z = \{ 0 \} $.

\item [(ii)] There exists an infinite subgroup $A \leq \T$ such that $
c_n \alpha  \to 0$ holds for all $ \alpha \in A$.
\end{itemize}
Conditions (i) is obviously equivalent to (i') below:
\begin{itemize}
\item [(i')] There exists a group compactification $(\iota,G)$
of $\Z$
such that $\iota $ is $1$-$1$,  $\iota (c_n) \to h$ in $G$
and $\langle h \rangle \cap \iota(\Z) = \{ 0 \}$.
\end{itemize}

We remark first that (ii) implies (i'):
Let $\iota:\Z\rightarrow \T^A, n\mapsto (n\alpha)_{\alpha\in A}$
and put $G=\overline {\iota(\Z)}$.
Obviously $(\iota,G)$ is a compactification of $\Z$ and since
$A$ is infinite, $\iota$ is $1$-$1$.
\\

To see that the converse does not hold, pick
$\alpha,\beta\in \T$, such that $\alpha$ and $\beta$ are linearly
independent over the rationals. Define a homomorphism
 $f:\langle \alpha \rangle
\rightarrow \langle \beta \rangle, n\alpha \mapsto n\beta$.
By Theorem \ref{appl} choose a sequence $\seq c$ in $\Z$
such that $c_n \alpha \rightarrow f(\alpha)=\beta$
and $(c_n \gamma)\io n $
does not converge for $\gamma \in \T\setminus \langle \alpha\rangle$.

Then $\iota:\Z\rightarrow \T, n\mapsto n\alpha$ gives rise to
a group compactification of $\Z$. Put $h:=\beta $
such that $ \iota (c_n)=c_n\alpha \rightarrow \beta$. Since
$\alpha$ and $\beta $ were chosen to be linearly independent
we have $\langle h \rangle \cap \iota(\Z) = \langle \beta\rangle
\cap \langle \alpha\rangle=\{0\}$. Thus $(i')$ is valid.
On the other hand $(ii)$ fails since $c_n \gamma\rightarrow 0 $
 only for $\gamma =0$.
\\

For a different type of counterexample
fix a prime $p$ and consider  the
     $p$-adic integers $\Z_p$.
Choose an arbitrary sequence $\seq k$ in ${\ohat {p-1}}$
 which contains  infinitely many non zero elements and
satisfies $k_1=1$. Using this, put for each $n\in \N$,
$h_{2n}=h_{2n+1}=\sum_{i=1}^n k_i p^i$ and
let $c_n = p^n+h_n$. Then
$$\limm n c_n = \limm n h_n = \sum_{i=1}^\infty k_i p^i=:h\in
\Z_p\setminus \Z.$$ Hence $kh\in \Z_p\setminus \Z$
for all $k\in \Z\setminus \{0\}$, so   $(i)$ holds.

Next pick $\alpha \in \T$ such that  $\limm n c_n \alpha=0$.
It follows that also $\limm n (c_{n+1} - c_n)\alpha=0$.
 Since $h_{2n} = h_{2n+1}$  this yields $p^n \alpha
\to 0$ so $\alpha= a/p^l +\Z$ for some $l \in \N$ and
 $a \in \ohat{p^l-1}$.
 But then, for all $n\geq l$
$$\| \alpha c_{2n}\| =
  \left\|(a/p^l) \left(p^{2n}+ \sum_{i=1}^n k_i p^i\right) \right\| =
\left\|(a/p^l) \sum_{i=1}^l k_i p^i \right\|.
 $$
The last term tends
to $0$ only if $a=0$. Hence $\alpha=0$ and $(ii)$ fails.

\end{document}